
\input vanilla.sty
\magnification 1200
\baselineskip 18pt
\input definiti.tex
\overfullrule=0pt
\input mathchar.tex
\define\pmf{\par\medpagebreak\flushpar}

\define\vp{\varphi}

\define\sii{\sigma_i}
\define\ai{\alpha_i}
\define\pbf{\par\bigpagebreak\flushpar}
\font\frak=eufm10

\define\a{\alpha}
\define\w{\omega}

\define\cala{{\cal A}}
\define\calb{{\cal B}}
\define\calf{{\cal F}}
\define\calc{{\cal C}}
\define\cald{{\cal D}}
\define\calh{{\cal H}}

\define\calp{{\cal P}}

\define\cals{{\cal S}}
\define\calw{{\cal W}}
\define\calz{{\cal Z}}
\define\calo{{\cal O}}

\define\qed{\hfill $\square$}

\define\C{ChS}

\centerline{{\bf Abstract}}

We prove the Bloch conjecture : $ c_2(E) \in H^4_\cald (X,\bbz(2))$ is torsion
for holomorphic rank two vector bundles $E$ with an integrable connection over
 a complex projective variety $X$.  We prove also the rationality of
the Chern-Simons invariant of compact arithmetic hyperbolic three-manifolds.
We give a sharp higher-dimensional Milnor inequality for the volume
regulator of all representations to $PSO(1,n)$ of fundamental groups of
compact $n$-dimensional hyperbolic manifolds, announced in our earlier paper.

\vfill\eject

\heading {\bf Rationality of Secondary Classes}
\endheading

\author {\bf Alexander Reznikov} \endauthor


\noindent{\bf 1. THE THEOREM}

\noindent{\bf 1.1.} Let $X$ be a smooth complex projective variety.
Consider a representation $\rho: \pi_1(X) \to SL(2,\bbc)$.  Let $E_\rho$
be the corresponding rank two vector bundle over $X$.
Viewing $E_\rho$ as an algebraic vector bundle, denote $c_2(E_\rho )$ the
second Chern class in Deligne cohomology group $H^4_\cald (X,\bbz(2))$ ([15],
[20]).Recall that there is an exact sequence $ 0\to J^2(X)\to H^4_\cald
(X,\bbz(2))\to H^4(X,\bbz(2))$ and  by the Chern--Weil theory, the image of
$c_2 (E_\rho )$ in $H^4(X,\bbz(2))$ is torsion. Therefore  $c_2(E_\rho )$ lies
in the image of $H^3(X,\bbc/\bbz)$ under the natural map $H^3(X,\bbc/\bbz)\to
H^3(X,\bbc/\bbz(2))\to H^4_\cald (X,\bbz(2))$. It was proved by Bloch [3],
Gillet--Soul\'e [24] and Soul\'e [49] that in fact, $c_2(E_\rho )$ is an image
of the secondary characteristic class $Ch(\rho)$ of a flat bundle $E_\rho $
(equivalently, of a representation $\rho$), lying in $H^3(X,\bbc/\bbz)$.The
$\bbr/\bbz$-part of this class was introduced and studied by Chern--Simons[9]
and Cheeger--Simons[8],and will be called Cheeger--Chern--Simons class and
denoted $ChS(\rho)$.The $\bbr$-part, lying in $H^3(X,\bbr)$ will be called
Borel hyperbolic volume class (regulator) and denoted $Vol(\rho)$. Remar!
k thast if $\rho$ is unitary, the
   Next, for a
field $F$ denote $\calb (F)$ the Bloch  group of $F$.  Recall that there is,
for $F$ algebraically closed, an exact sequence $0\to \mu^{\otimes 2}_F \to
H_3 (SL (2,F), \bbz) \to \calb(F) \to 0$ of Bloch--Wigner--Dupont--Sah [18].
 The dilogarithm function of Bloch--Wigner defines a homomorphism
$D: \calb (\bbc) \to \bbc/\bbq = \bbr/\bbq \oplus i \bbr$ which splits to
the Borel hyperbolic volume regulator and the Bloch--Beilinson Chern-Simons
regulator
maps [19].  It is proved by Bloch [3] that for $\rho$ unitary, the reduction
of $Ch(\rho)$  $mod\bbq$ lives in
$\frac{Span(Re D (\calb(\bar \bbq))}{\bbq}$, hence  $c_2(E_\rho)$ assumes at
most countably many
values in $H^4_\cald (X,\bbz(2))$.

In [3], Bloch went on to conjecture the following result, which will be proved
here.

\proclaim{\bf  Main Theorem A}(\bf The Bloch conjecture). For any
representation $\rho:\pi_1(X)\to SL(2,\bbc)$,  the class $c_2(E_\rho )\in
H^4_\cald (X,\bbz(2))$ is torsion.
\endproclaim

The proof of the Main Theorem A will be completed in the section 4.
The strategy we choose is described in the following lines.

\noindent{\bf 1.2.} The main ingredients of the proof are: the rational
algebraic
$K$-theory, the homological finiteness of $S$-arithmetic groups, the
existence theory of twisted harmonic maps, and the Siu--Sampson--Carlson--
Toledo rigidity theory.

First, we use the above-cited work of Bloch [3], Gillet--Sou\'le [24] and
Soul\'e [49]
to reduce the theorem to the following two statements: for any $\rho$ as above
i)$Vol(\rho)=0$ and ii)$ChS(\rho)\in H^3(X,\bbq/\bbz)$.
Next, we consider  the  representation variety  $Hom (\pi_1 (X), \;
SL_2)= V_X$.  This is a scheme over $Spec (\bbz)$.  Let $V_X(p)$ be the
irreducible component of $V_X(\bbc)$, containing $\rho$.  Then
$V_X(\rho)$ contains a $\bar \bbq$-point, say $\bar \rho$.  Since $V_X(\rho)$
is
connected in the classical topology, the rigidity of the Chern--Simons  and
Borel classes
[9] gives $Vol(\bar \rho)=Vol(\rho)$,  $Ch S(\bar \rho)=Ch S(\rho)$. So we may
assume that $\rho$ is
defined over an $S$-arithmetic groups $\calo_S \subset F$, where $F$ is a
number field.  Consider the induced map $\hat \rho: X \to BSL_2 (\calo_S)$.  We
will show, using the work of Borel and Serre [6] on the finiteness
properties of $S$-arithmetic groups, that there exists a universal class $C:
H_3 (SL_n (\calo_S))\to \bbr $, $n\geq2$, and a natural number $M$ such that
$M\cdot ChS(\rho)=\hat \rho^*(M\cdot C)(mod \bbz)$.

Let $\sigma_1, \dots \sigma_m$ be the maximal set of nonconjugate embeddings
of $F$ into $\bbc$.  Let $Vol\in H^3 (BSL(\bbc),\; \bbr)$ be the universal
Borel
hyperbolic volume regulator.  Then by [4], the group $H^3 (BSL_n(\calo_S),
\bbr)$ is freely generated by $\sigma^*_i Vol, \; i=1, \dots m$,for $n\geq12$.
We will
prove, using the theory of harmonic sections and a version of the degeneration
result of Sampson, that $Vol (\sigma_i \circ \rho)$ vanishes for all $i$.
Hence $\hat \rho_* (H_3 (X, \bbz)) \subset H_3 (SL_n (\calo_S))$ is the torsion
part of $H_3 (SL_n (\calo_S))$, so $M\cdot ChS(\rho)$ is zero in
$\bbc/\bbz$, and $ChS(\rho)\in H^3(X, \bbq/\bbz)$.

Using the different geometrical argument (the theory of Gromov' simplicial
volume invariant) we will prove in the section  5 the  following sharp higher
Milnor
inequality for volume invariants,announced in [41].
\proclaim {\bf Theorem B} Let $M$ be a compact $n$-dimensional hyperbolic
manifold and let $\mu:\pi_1(M)\to PSO(1,n)$ be a representation. Then
$$Vol(\mu) \leq Vol(M)$$
\endproclaim

   We will also deduce the following result
whose evidence was based on the computations of Fintushel and Stern [21], and
Kirk and Klassen [32]:
\proclaim{\bf Theorem C} Let $M^3$ be a Seifert fibration, and let $\rho: \pi_1
(M) \to SL_2 (\bbc)$ be a representation.  Then $(Ch S(\rho),[M])\in \bbr/\bbz$
is rational.
\endproclaim

Finally, we will prove the following result.
\proclaim{\bf Theorem D} Let $M^3$ be compact  arithmetic hyperbolic manifold
and
let $\rho: \pi_1 (M) \to PSL_2 (\bbc)$ be the defining representation.
Then $(ChS(\rho),[M])$ is rational.
\endproclaim

I am deeply grateful to Professors John Jones, Mel Rothenberg, Ron Livne,
Hel\'en\`e Esnault, Don Zagier, Norbert Schappacher, Udi de Shalit
and Shicheng Wang for very helpful discussions on the subject of this paper.
I am also grateful to Professors Jim Eells and J\"urgen Jost  for numerous
discussions on the harmonic maps theory.Special thanks are due to H.Esnault and
the referee for numerous improvements and critisism of various inaccuracies in
the first version of the paper.



\noindent{\bf 2. REDUCTION TO $S$-ARITHMETIC GROUPS}

\noindent{\bf 2.1.}Let $X/\bbc$ be as in 1.1. and let $\rho:\pi_1(X) \to
SL_n (\bbc)$ be a representation.
The following result was proved by
Bloch [3], Gillet-Soul\'e [24], and Soul\'e [49]:

\proclaim{Theorem} The Chern class $c_i(E_\rho )\in  H^{2i}_\cald (X,\bbz(i))$
is the image of the secondary class $Ch(\rho )\in H^{2i-1}(X,\bbc/\bbz)$ under
the natural map $H^{2i-1}(X,\bbc/\bbz)\to H^{2i-1}(X,\bbc/\bbz(i))\to
H^{2i}_\cald (X,\bbz(i))$.
\endproclaim

\noindent{\bf 2.2.} For any finitely-generated group $\Gamma$ and any
algebraic group $G/\bbq$ let $V^G_\Gamma$ be the representation variety
$Hom(\Gamma, G)$.  This is an affine scheme defined over $Spec(\bbq)$ [12],
[29].
Hence any irreducible component of $V_\Gamma (\bbc)$ contains a $\bar
\bbq$-point.For  $G=GL_n$ and $\rho:\Gamma \to GL_n(\bbc)$ let $V_\Gamma
(\rho)$ be the component, containing $\rho$ and let $\bar \rho$ be a $\bar
\bbq$-point in $V_\Gamma(\rho)$.By the rigidity of the secondary classes [9] we
have  $ChS(\rho)=ChS(\bar \rho)$ and $Vol(\rho)=Vol(\bar \rho)$.So for proving
the Main Theorem
we may assume that $\rho$ is defined over $\bar \bbq$ and we need to show
that for any singular manifold $i:M^3 \to X, \;(Vol(\rho \circ i_*),[M])=0$ and
$(ChS(\rho \circ i_*),[M])$ is rational.
Since $\pi_1(X)$ is finitely-generated, $\bar \rho$ is actually defined over
$\calo_S$, where $S$ is a finite set of primes.

\noindent{\bf 2.3.} Let $BSL_n(\calo_S)$ be the classifying space of
$SL_n(\calo_S)$.  A representation $\rho:\pi_1(X)\to SL_n(\calo_S)$ determines
a homotopy class of maps $\hat \rho:X \to BSL_n(\calo_S)$, and conversely.  By
the deep theorem
of Borel and Serre [6], $H_*(BSL_n(\calo_S)$ is finitely generated.  In the
following section we will develop for flat bundles over such spaces a general
theory of regulators, and show the existence of the universal $\C$-class
in $H^3(BSL_n(\calo_S),\bbr/\bbz)$.

Observe that one can define a universal class in
$H^3(BSL^\delta_n (\bbc),\bbr/\bbz)$ following Beilinson and Soul\'e, by using
the natural map of simplicial schemes $BSL^\delta_n(\bbc)_.\to
BSL_n(\bbc)_.$.Indeed, $H^{2n}(|BSL_n(\bbc)_.|)$ has pure $(n,n)$-type [15] and
so the usual Chern class of the classifying bundle lifts to Deligne cohomology.
 However, we need a by-hand construction of 3.3.
below, since it fits the geometric framework, in which various regulators
are constructed in 3.2.1.--3.2.4.  This paves the way for our use of heavy
analytical weapons in the section 4, and ultimately is one of the most
important components of the success in proving the Bloch conjecture.
The interpretation of the Chern-Simons regulator and the
Borel hyperbolic volume regulator given (within the scope of a more general
theory) below in the section 3, was largely inspired by the
theory of characteristic classes of foliation, developped by Bott-Haefliger
and Bernstein-Rosenfeld in parallel with the ground-breaking work of Chern
and Simons.


\noindent {\bf 3. GENERAL THEORY OF REGULATORS}

\noindent{\bf 3.1.} Let $\calc\calw$ denote the category of CW-complexes and
let $\calh:
\calc\calw\to \cala b$ be a cohomology theory.  A functor $\calz:
\calc\calw \to \cals ets$ and a
morphism  [\ ]: $\calz\to \calh$ form a cocycle theory for $\calh$.  If for any
$X
\in \calc\calw$ the presheaf $U\mapsto \calz(U)$ on $X$ is a sheaf, we say
that the theory $(\calz, [\;])$ is infinitesimal.  For example, the theory
$X\mapsto$ (maps from $X$ to $K(\pi ,1)) $ is infinitisemal,and $X\mapsto$
(singular cocycles on $X$) is not, both attached to the singular
cohomology theory.  For a subcategory $C^\infty$ of smooth manifolds we have
a cocycle theory $\Omega:X \mapsto$ (closed exterior forms on $X$), which is
infinitesimal.

Next, let $\pi_1(Y) \buildrel \rho \over \to Homeo (X)$ be a representation,
and let $z \in \calz(X)$ be a cocycle, invariant under $\rho(\pi_1 (Y))$.
Consider the flat bundle $X\to F_\rho \to Y$, where $F_\rho=\buildrel \sim
\over  Y \mathop
\times \limits_{\pi_1 (Y)} X$.  If $Y$ is locally simply-connected, this bundle
is
locally trivial.  Let $\mathop \bigcup\limits_i U_i=X$ be a covering of $X$ by
opens, such that $\pi_1(U_i)=0$, so that $F_\rho|U_i$ is trivial.  Fix an
identification $F_\rho|U_i {\mathop{\tilde\to}\limits^{\a_i}} U_i \times X$
coming
from the flat connection.  Denote $p_2$: $U_i \times X\to X$ the projection to
the second
factor. This induces an element $(\a_i \circ p_2)^* z \in \calz
(F_\rho|U_i)$, denoted $y_i$.
Since $z$ is invariant under $\rho (\pi_1(Y))$, these $y_i$ form a compatible
family, so if $\calz$ infinitesimal, there is a well-defined element $y\in
\calz(F_\rho)$.

Assume that $X$ is contractible.  Then all sections $s:Y\to F_\rho$
are homotopic, and we obtain a well-defined element $[s^*y]\in \calh(X)$ called
the regulator of $\rho$ and denoted $r(z, \rho)$.  If $X$ is not
contractible, then anyway we obtain an element $[y]\in \calh(F_\rho)$ and then
use the spectral sequence of the fibration $X\to E_\rho\to Y $ to
get secondary invariants of $\rho$.

\noindent{\bf 3.2. Examples}

\noindent{\bf 3.2.1. Classical Borel Regulators}  Let $G$ be a real Lie group
and let $K\subset G$ be a maximal compact subgroup of $G$.  For any manifold
 $Y$ and a representation $\rho: \pi_1(Y) \to G$ we get a homomorphism
$\wedge^*_K (\hbox{\frak\$g/\frak\$k})\buildrel Bor \over \to H^* (Y)$,
described as follows.
Let $x \in \wedge^*_K(\hbox{\frak\$g/\frak\$k})$.  Define a $G$-invariant
form on $G/K$,
corresponding to $x$ and use the construction of 3.1.  In particular, if
$(G,K)=(SL_n(\bbc), SU(n))$, then the image of the properly normalized
generator of $\wedge^3_K(\hbox{\frak\$s\frak\$l}_n/\hbox{\frak\$s\frak\$u}(n))$
is called the
hyperbolic volume of
$\rho$, denoted $Vol_3(\rho)$.  For the pair $SO(1,n), \; SO(n))$ we get
an element $Vol_n(\rho)$ corresponding to the generator of $\wedge^n
(\hbox{\frak\$s\frak\$o} (1,n)/\hbox{\frak\$s\frak\$o}(n))$, also called the
hyperbolic volume, see [41], for example.

\noindent{\bf 3.2.2. Cheeger--Chern-Simons Classes}  Again let $G$ be a real
Lie
group, and let $x \in \wedge^* (\hbox{\frak\$g})$.
Let $z$ be a left-invariant form in $\Omega^*(G)$, corresponding to $G$.
Then for any manifold $Y$ and any
representation we get an element $[y]\in H^* (F_\rho)$.  If $G$ is
contractible,  e.g. $G=\tilde{SL_2} (\bbr)$, this defines an element
$r(x, \rho)\in H^* (Y,\bbr)$ as above.  We refer to [44], to
the detailed study of this last example.  If $G$ is not contractible, one
looks at the spectral sequence of $F_\rho$ to see what can be done to descend
some cohomology information down to $Y$.

Let us specialize this construction for the Cartan form $\w(X,Y.Z)=([X,Y],Z)$,
where $(\cdot,\cdot)$ is the Cartan-Killing scalar product in $sl_n(\bbc)$.
This is a complex-valued invariant 3-form, so that we may look at $r (Re\; \w,
\rho)$ and $r (Jm \;\w, \rho)$ both in $H^3 (F_\rho, \bbr)$.  Observe that
$F_\rho$ here is a flat principal $SL_n(\bbc)$-bundle over $Y$.  Now, since
$Jm\;\w$ is exact, we may descent $r(Im \;\w,\rho)$ to $H^3(Y,\bbr)$.We claim
 that this will be precisely the hyperbolic volume
regulator of 3.2.1.Indeed, fix a point $p\in \calh^3$and consider the
evaluation map from $SL(2,\bbc )$ to $\calh^3$, sending $g$ to $gp$.This map is
equivariant with respect to $SL(2,\bbc )$- actions considered and the pull-back
of the volume form on $\calh^3$ is precisely $Jm \;\w$, and we can used the
functoriality by $X$ in 3.1.  On the other hand, $Re\;\w$ is not exact and
represents a generator of $H^3(SL_n(\bbc),\bbr)\approx \bbr$.  Normalizing it
to $\frac{1}{4\pi^2}Re\;\w$ so that its period will be one, we easily see, in
the case when $F_\rho$
is topologically trivial, that the descend by different sections will give a
well-defined element in $H^3(Y,\bbr/\bbz)$.  This is the classical
Chern-Simons  class.If $F_\rho$ is not topologically trivial, then it is  still
well-defined [9].
\pmf
{\bf 3.2.3. Thurston-MacDuff-Morita invariants.}
Let $X$ be a contractible manifold with a volume form $\w$, e.g.$ (X,\w)=
(\bbr^n, can)$.  Then for any manifold $Y$ and a representation $\rho: \pi_1
(Y) \to ${\it Diff}$_{\w}(X)$ one applies 3.1 and gets an element
$Bor (\w, \rho)\in H^n(Y,\bbr)$.
It is easily shown that this element comes from the
universal class $Vol_\w \in H^n$ ({\it Diff}$_{\w}(X),\bbr)$.  Similarly let
$(X, \sigma)$ be a contractible symplectic manifold, then there exists and
element {\it Sympl} $ \in H^2(B${\it Sympl} $ (X), \; \bbr)$, where
{\it Sympl} $(X)$ is the symplectomorphism group.  These classes were defined
in
[34], [35], [53], [31], using simplicial constructions.
A somewhat deeper look at the topology of symplectic fibrations enables one
to define a Chern-Simons-type invariants, c.f. [46].
\pmf
{\bf 3.2.4. $K$-Theoretic Invariants of Group Actions}  Let $\Gamma$ be a
group acting smoothly on a manifold $X$.  Form a flat bundle $X\to F \to
B\Gamma$ and consider a vector bundle $\calf$ over $F$, tangent to fibers.
This is well-defined despite the fact that $B \Gamma$ is not a manifold.
If $X$ is contractible, find a section $s$ of $F$ and consider the class
$s^*[\calf]\in K^0(B \Gamma)$, which is a well-defined invariant of the
action.  If $X$ is not contractible we may look at the characteristic class of
$\calf$ in the singular cohomology of $F$ and try to get secondary invariants
in $H^*(\Gamma)$, using the spectral sequence, as above.  This construction
does not fall under the axiomatic description of 3.1. since the cocycle theory
$X\mapsto$ (isomorphic classes of vector bundles over $X$) is not
infinitesimal.  The problem of finding better axiomatic description, covering
this case, is left for the interested reader.

The secondary classes of group actions, described above, prove important in
applications to the finite group actions.  Details will appear elsewhere.
\pmf
{\bf 3.3.} In this paragraph we will describe the construction of regulators
for the de Rham cocycle theory in the case when $Y$ is a locally simply
connected $CW$--complex, not a
manifold (comp.[19]). We concentrate on the special case $X=SL_n(\bbc)$.  It is
always
assumed that $H_*(Y) $ is of finite type.  Then the image of the natural map
$MSO_k(Y)\to H_k(Y)$ is of finite index in $H_k(Y)$ by the theorem of Thom.
For any representation $\rho:\pi_1(Y)\to SL_n(\bbc)$ the usual Chern-Weil
theory will imply then that all Chern classes $c_i(E_\rho) \in H^{2i}
(Y, \bbz)$ die after tensoring by $\bbr$, hence $c_i(E_\rho) \in
H^{2i}_{tors} (Y, \bbz)$.Let $N$ be such that $Nc_2(E\rho )=0$, then
$c_2(E_\rho)$ is a Bockstein image of a class $z$ in $H^3(Y, \bbz/N\bbz)$.Let
$f:Y\to K(3,\bbz/N\bbz)$ be the classifying map and let $\bar Y$ be a
homotopical fiber of $f$.Then the inclusion $i:\bar Y \to Y$ induces an
isomorphism in the rational homology.  Pull back the flat bundle $E_\rho$ from
$Y$ to $\bar Y$and denote $\bar E_\rho$ the resulting bundle. Then $c_2(\bar
E_\rho)=0$, hence the principal $SL_n(\bbc)$-bundle $\bar F_\rho$, assosiated
to $\bar E_\rho$ admits a section over the 4-skeleton $Sk_4(\bar Y)$.  Fix
such a section $s$.  Let $\bar \rho= \rho \circ i_*: \pi_1(\bar Y) \to SL_n
(\bbc)$, then $\bar E_\rho$ is the flat bundle, associated with $\bar \rho$.
Consider a singular manifold $j:M^3\to{Sk_4\bar  Y}$ and a flat bundle
$j^*\bar E_\rho$ with the canonical {\it smooth} structure (of a flat
bundle).  Consider the 3-form $x=Re([X,Y],Z)$ on $sl_n(\bbc)$, where
$(\cdot ,\cdot)$
is the (complex) Cartan-Killing scalar product.  The corresponding
regulator $r (x, j^*\rho) \in H^3 (j^*\bar F_\rho, \bbr)$ is the
usual Chern-Simons invariant.  Next, we use  {\it the canonical} choice of
a section, namely, $s\circ j$, to produce a class, also called $r
(x, j^*\rho)$, in $H^3(M,\bbr)$, and a number $(r (x, j^*\rho),
\;[M])\in \bbr$.  We claim this defines a homomorphism  $MSO_3(\bar Y) \to
\bbr$.
Indeed, suppose we are given a map $\psi:N^4\to \bar Y$ with $(M,j)=\partial
(N,\psi)$.
Arguing as above, we get a class $r (x, \psi^* \rho)\in H^3 (\psi^*
\bar F_\rho,\bbr)$, whose restriction on $\phi^*\bar F_\rho|_{\partial N}$
gives $r (x, j^*\rho)$.  Then it is obvious that the latter class is zero.

Since $H_* (\bar Y)$ is finitely generated we get a homomorphism $H^3
(\bar Y)\approx H_3 (Sk_4\bar Y)\to \bbr$, whose reduction $mod\; \bbz$
induces a
usual Chern-Simons class on every singular manifold $j: M^3 \to \bar Y$.
Now, since $i_*:H_3(\bar Y,\bbq)\buildrel \sim\over \to H_3(Y,\bbq)$,
we get a homomorphism $C:H_3(Y,\bbz)\to \bbr$, such that for some number
$M\in \bbn$ big enough the following is true:  for any singular manifold
$j:M^3 \to Y,\quad C(j_*[Y])=ChS(j\circ \rho) \quad (mod \bbz \cdot
\frac{1}{M})$.
\pmf
{\bf 3.4. Remark.} As it was mentioned above,there exists a universal $ChS$
class in $H^3(BSL_2(\bbc),\bbr/\bbz))$.



\noindent{\bf 4. PROOF OF THE MAIN THEOREM}

\pbf {\bf 4.1}  Let $F$ be a number field without real places and let
$\sigma_1,\dots
\sigma_m$ is a maximal family of nonconjugate embeddings of $F$ into $\bbc$.
Let $\calo_S$ be as above.  Consider the universal Borel regulator $Vol\in H^3
(BSL_n(\bbc), \bbr)$.  We need the following fundamental result of Borel [4].
\proclaim{Theorem} (A. Borel) The elements $\sigma_i^*\hbox{Vol}\;
\in H^3(BSL_n (\calo_S),
\bbr)$ form a basis of $H^3(BSL_n (\calo_S),\bbr)$ over $\bbr$ for $n$
big enough $(n\ge 12)$.
\endproclaim
\pmf {\bf 4.2.} Combining the theorem 4.1. with 3.3., we get the following:

\proclaim{Fundamental lemma} There exist constants $\a_1,\dots \a_m \in \bbr$
and a natural number $M'$, such that for any representation $\rho:\pi_1(M^3)\to
S L_n
(\calo_S)$ one gets
$$
ChS(\rho)\equiv\sum\limits^m_{i=1} \a_i Vol (\sigma_i \circ \rho)\quad (mod
\bbz\cdot \frac{1}{M'})\tag *
$$
where $n$ is as in 4.1.

Since the $ChS$-invariant is compatible with the embeddings $SL_m\to SL_n,\;\;
n>m$, we can remove the restriction on $n$.
\endproclaim
\pbf {\bf 4.3.} To prove theorem A we need to show that for any
representation $\mu: \pi_1(X)\to SL_2 (\bbc),\quad Vol (\mu)=0$.  Applying
this to $\mu=\sigma_i\circ \rho$ we will get $ChS(\rho)\in  \bbz\cdot
\frac{1}{M'}$ on $MSO_3(X)$ by 4.2.

Consider the natural action of $SL_2(\bbc)$ on the hyperbolic space $\calh^3$.
Let $S_\infty$ be the sphere at infinity.  We may assume that $\mu(\pi_1(X))$
does not have fixed points in $S_\infty$, otherwise the
representation $\mu$ factors through $\bbr_+\times (SO(2)\ltimes \bbr^2)$ and
so deformes to a representation $\bar \mu$ in $SO(2)\ltimes \bbr^2$ and
$Vol(\mu)=Vol(\bar \mu)=0$ by
4.5.1 and rigidity of $Vol$. Consider the flat bundle $\calh^3\to \calf \to X$
associated to $\mu$.
We need the following fundamental result of Donaldson and Corlette [17], [13],
see also [33].

\proclaim{Theorem} (S. Donaldson, K. Corlette) Let $X$ be a compact Riemannian
manifold and let $\mu:\pi_1(X)\to SL_2(\bbc)$ be a representation.  If the
action of $\pi_1(X)$ on $S_\infty$ is fixed-point-free, then the flat bundle
$\calf$ posesses a harmonic section.
\endproclaim
\pbf
{\bf 4.4.} Harmonic maps of K\"ahler manifolds to manifolds of negative
curvature were intensively studied by many authors, starting with the seminal
work of Siu [48].  We quote in particular the following remarkable result of
Sampson [47].
\proclaim{Theorem (Sampson)}  Let $X$ be a compact K\"ahler manifold and let
$Z$ be a (real)
hyperbolic manifold.  Let $\phi: X \to Z$ be harmonic, then rank $D\phi \le 2$
everywhere on $X$.
\endproclaim

The proof of this theorem remains valid for harmonic sections of flat
bundles, associated with a representation $\mu: \pi_1(X)\to Iso (Z)$, as
in [41].Indeed, the proof in [47] is based on the Bochner--type integral
formula, which is local in  $\phi$ and hence holds for harmonic sections.
\pbf{\bf 4.5.}  Now we are ready to complete the proof of the theorem A.
Let $\mu:\pi_1 (X) \to SL_2(\bbc)$ be a representation.  Assuming that
$\pi_1(X)$ acts fixed-point-free on $S_\infty (\calh^3)$, consider a harmonic
section $s$ of $\calf$ which exists by 4.3. Let $\w \in \Omega^3(\calf)$ be
the three-form defined by the pull-back from the volume form of hyperbolic
fibers; it is well-defined by 3.1. The form $s^*\w$ represents $r (\mu)=
Vol(\mu)\in H^3 (X,\bbr)$. On the other hand, $Ds_x$, viewed as a map from
$T_xX$ to $T_{s(x)} \calf_x$ has a rank $\le 2$ by 4.4., for all $x \in X$,
so $s^*\w$ vanishes identically.  This proves $Vol (\mu)=0$ and hence
$ChS(\rho)\in H^3(X,\bbz\cdot \frac{1}{M'}/\bbz)$ by 4.2.
\pbf {\bf 4.5.1.}  So what is left is to show that if the representation $\mu
: \pi_1(X)\to SL_2 (\bbc)$ factors through $SO(2) \ltimes \bbr^2$, then $Vol
(\mu)=0$.  Recall [40] that the volume invariant may be interpreted as follows.
Consider the continuous cohomology $H^*_c(SL_2(\bbc))$.  There exists a
natural map $H^3_c(SL_2(\bbc)_,\bbc)\to H^3 (SL_2^\delta (\bbc),\bbc)$.  The
left hand side space has dimension 1, and the image of the canonical generator
is
precisely the element $Vol \in H^3 (SL_2^\delta (\bbc))$.  Now, we have a
diagram
$$
\matrix H^3_c (SL_2(\bbc), \bbc) &\longrightarrow &H^3_c (SO(2)\ltimes
\bbr^2)^\delta,\bbc)\\
\downarrow &&\downarrow\\
H^3(SL_2^\delta(\bbc),\bbc) &\longrightarrow &H^3 ((SO(2)\ltimes
\bbr^2)^\delta,\bbc)
\endmatrix
$$
Next, since $SO(2)$ is compact and $\bbr^2$ is acyclic we have the canonical
isomorphism [21] $H^3_c (SO(2)\ltimes \bbr^2,\bbc)\simeq H^3
(\hbox{\frak\$s\frak\$o}(2)\ltimes
\bbr^2, SO(2))=0$.
\pbf {\bf 4.6.} Assume that $X$ is a compact K\"ahler manifold and
$\rho:\pi_1(X) \to SL_2(\bbc)$
is a representation.  Then $ChS(\rho)$ is rational.   The proof is identical
to  that of theorem A.  Alternatively, for discrete representations one can use
the theorem A and the
result of Mok [38] which says that any such representation $\rho$ factors
through a homomorphism $\pi_1(X)\to \pi_1(Z)$ with $Z$ smooth and
projective.  On the other hand, the theorem A probably fails for compact
complex manifolds $X$ in view of the recent result of Taubs, saying that any
finitely presented group is a complex manifold group.

The result on vanishing of $Vol (\mu)$, where $\mu:\pi_1(X)\to SL_2(\bbc)$
is not true for representations in Lie groups of higher rank.
However, there are important results on rigidity of representations with
prescribed higher dimensional volume [14].  This also uses the harmonic
section technique.  Another interesting application is a construction of a huge
family of compact symplectic manifolds, which do not admit any K\"ahler
structure and have a fundamental group of an exponential growth.  See [43] for
the details.

\pbf{\bf 4.7.}  The main theorem above provokes a natural question:  what
happens for flat bundles of higher rank?  It is clear, because of what we
have just observed, that the proof presented here does not apply directly
to this situation.  Still, one can expect the validity of the following
statement.

\noindent{\bf  Bloch conjecture, higher rank case}  Let $X$ be a complex
projective variety
and let $\rho: \pi_1 (X) \to SL_n (\bbc), \;\;n\ge 2$, be a representation.
Let $E_\rho$ be the corresponding flat rank $n$ vector bundle.Then for all
$i\geq 2$, $c_i (E_\rho )\in H^{2i}_\cald (X,\bbz(i))$ is torsion.

As stated, this conjecture becomes very similar to the generalized
Bloch-Beilinson conjecture on higher regulators of motives.  Anyway, this
conjecture is true; the detailed proof and all the yoga around it will be
presented in the forthcoming paper [45].

On the other hand, for $i=1$ and the representations in $GL_n(\bbc)$ one
cannot expect anything like this to hold, just because even for curves, the
first Chern-Simons class of a flat line bundle is just its monodromy
presentation $\pi_1(X)\to \bbc^*$, which may be completely arbitrary.
However, if $X$ is defined over a number field, $k$ and if the flat
bundle comes from the Deligne-Ramakrishnan construction, associated to
an element $z$ of $K_2(X_k)$ [40], one expects a lot of rigidity for
the monodromy representation.  The celebrated Bloch-Beilinson conjecture
relates its periods to the value of $L$-function of $X$ at $s=0$.
Moreover, the logarithms of the real parts of these periods are believed to
form a $\bbq$-structure of $H^1(X,\bbr)$ if $[k:\bbq]=1$ and $z$ varies
inside $K_2(X_\bbq)\otimes \bbq$.  The reader will find the results of our
geometrical treatment of higher regulators in algebraic $K$--theory in a
forhcoming paper.


\noindent{\bf 5. HIGHER MILNOR INEQUALITY AND $ChS$ INVARIANTS OF
SEIFERT FIBRATIONS}

\noindent{\bf 5.1.} In this section we will prove the theorems B and C.  The
proof will be
based on the sharp higher Milnor inequality, announced in [41] with the
proof of a weaker estimate given there.  For a manifold $M$ we denote
$\| M\|_g$ the Gromov's simplicial volume of the fundamental cycle.
Let $\mu: \pi_1(M) \to PSO(1,n)$ be a representation. Since
$PSO(1,n)$ acts isometrically in $\calh^n$,  the general theory of 3.1
gives us an invariant $Vol(\mu)\in H^n(M,\bbr)$. For $n$-dimensional $M$ we
denote $(Vol (\mu), [M])$ again by $Vol (\mu)$.  Then we
state the following result:
\proclaim{Theorem} For any compact $n$-dimensional manifold $M$ and any
representation $\mu: \pi_1(M) \to PSO(1,n)$ the volume $Vol(\mu)$
satisfies
$$
Vol(\mu) \le d_n \|M\|_g,
$$
where $d_n$ is the Milnor constant, i.e., the volume of the regular infinite
simplex.
\endproclaim
\pbf {\bf 5.2.} Combining 5.1. with a result of Gromov [27], [52] we get
\proclaim{Theorem B} (Higher Milnor Inequality).  Let $M^n, n\ge 2$ be a closed
hyperbolic manifold and let $\mu: \pi_1(X)\to PSO(1,n)$ be a representation.
Then
$$
Vol(\mu)\le Vol\; M.
$$
\endproclaim

The discussion of the classical case $n=2$ with various generalizations is
to be found in [41].  The proof of 5.1., 5.2. follows the pattern given in
[41].
\pbf{\bf 5.3. Remark.} The theorems 5.1., 5.2. are directly inspired by, and
are generalizations of the inequality of Gromov-Thurston [52] saying that for
any map $\phi: M \to N$ with $N$ hyperbolic one has
$$
deg \; \vp\cdot Vol\; N \le d_n\cdot \| M \|_g
$$
\pbf {\bf 5.4.} Using 5.1. we now prove the theorem C.  Let $M^3$ be a
Sifert fibration by a theorem of Yano [56], $\|M\|_g=0$.  Hence by 5.1.
$Vol(\mu) =0$ for any $\mu$.  Then applying 4.2. we get $ChS(\rho) \in H^3
(M, \bbz \cdot \frac{1}{M'}) $. \qed
\pbf{\bf 5.5.} \proclaim{Example} Let $\Gamma \subset \tilde{SL_2}
(\bbr)$ be a uniform lattice.  Then for some canonical choice of the Haar
measure in $\tilde{ SL_2}(\bbr)$,
$$
Vol (\tilde{SL_2} (\bbr) /\Gamma) \in \bbq
$$
\endproclaim
\demo{Proof}  We normalize the $ChS$ form on $SL_2(\bbc)$ in such a way
that the period of this form (the integral over $SU_2 \subset SL_2 (\bbc))$ is
one.  For a uniform lattice $\Gamma \subset \tilde{SL_2} (\bbr)$ denote
$M=\tilde{SL_2} (\bbr)/\Gamma$ and let $\rho: \pi_1(M) =\Gamma \to
\tilde{SL_2} (\bbr) \to SL_2 (\bbr) \to SL_2(\bbc)$ be the canonical
representation.
Then Covol ($\Gamma) =ChS(\rho)$.  Moreover, $M$ is a Seifert fibration
[37] and 5.4. applies to prove 5.5.
\pbf {\bf 5.6. Example.} Let $M=\{ z_1^p +z_2^q +z^r_3 =0 \} \cap S^{5} \subset
\bbc^3$ be a Pham-Brieskorn homology sphere ($p,q,r$ are coprime and $\frac
{1}{p} + \frac{1}{q} +\frac{1}{r} < 1$). Then by Milnor [37] and Dolgachev
[10],
$M=\tilde{SL_2}(\bbr)/\Gamma$, so that $Vol(M)$ is rational.
The $ChS$-invariant of representations in $SL_2(\bbr)$ has hyperbolicity
properties studied in [44].  In particular, it is shown there that the
topological complexity of $\mathop\#\limits_m M$ grows linearly with $M$
if $M$ is a Pham-Brieskorn variety.
\pbf {\bf 5.7. Example.} Let $M=\sum (a_1, \dots a_n)$ be a genus zero Sifert
fibration, so a homology sphere.  The representation variety $Hom (\pi_1 (M),
SU_2)/SU_2$ is studied by Fintushel-Stern[21], Kirk-Klassen [32],
and Bauer-Okonek [1].  Fintushel and Stern [21] computed $ChS$ invariants of
all unitary representations of $\pi_1(M)$.  These are rational numbers with
denominator a divisor of $a=a_1 \dots a_n$.  It would be interesting to
compute $ChS$ invariants of all representations in $SL_2(\bbc)$.
\pbf{\bf 5.8.} Let $M$ be a compact oriented three-manifold.  We will say that
$M$ is {\bf of hyperbolic type}, if there exists a representation $\mu:
\pi_1(M)\to SL_2 (\bbc)$ with $Vol(\mu) \neq 0$.  We actually have proved the
following:
\proclaim{Proposition} If there exists a representation $\rho: \pi_1(M) \to
SL_2 (\bbc)$ with irrational $ChS$-invariant, then $M$ is of hyperbolic type.
\endproclaim

Moreover, as we just have seen, all Seifert fibrations are not of hyperbolic
type.  We notice that the manifolds of hyperbolic type enjoy many properties
of actually hyperbolic manifolds, as follows.
\pbf {\bf 5.9.}\proclaim{Proposition} (comp. [27], [52]). Let $N,M$ be
compact oriented three-manifolds such that $M$ is of hyperbolic type.  Then
the degrees of continuous maps $\phi: N\to M$ are bounded by a constant
$C(N,M)$.
\endproclaim
\demo{Proof}  Let $\mu: \pi_1(M)\to SL_2 (\bbc)$ be a representation with
$Vol(\mu) \neq 0$.  For a map $f: N\to M$ we have by the naturality of the
volume invariant: $Vol (\mu \circ f_*)=\deg f \cdot Vol (\mu)$.  On the
other hand, $Vol (\mu \circ f_*)\le \| [N]\|_g$ by 5.1. Hence $\deg f \le
| Vol (\mu)|^{-1} \| [N]\|_g$.
\pbf{\bf 5.10.} \proclaim{Corollary} (c.f. [27],[52]).  Let $M$ be of
hyperbolic
type.  Then for a self-map $f: M \to M$, $\deg f$ assumes one of the values:
$0, \pm 1$.
\endproclaim
\pbf{\bf 5.11. Remarks}(i) If $M$ is of hyperbolic type, and $N$ is
any three-manifold, then $M \# N$ is of hyperbolic type
\item\item{(ii} If $M$ is of hyperbolic type, and $\vp: N \to M$ is a ramified
covering along a link in $M$, then $N$ is of hyperbolic type.

Combining the operations (i) and (ii) one constructs in abundance
{\bf non-hyperbolic irreducible} manifolds $M$, for which the analogs of the
Gromov-Thurston theorems (5.9. and 5.10.) are still true.
\pbf {\bf 5.12.} A well-known theorem of Goldman [25]  asserts that for $n=2$,
the equality in 5.2. takes place iff the monodromy group of $\mu$ is a uniform
lattice in $PSL_2(\bbr)$.  This generalizes to $n \ge 3$ as follows: any
representation for which $Vol (\mu)= Vol (M)$ is conjugate to a
composition $\rho \circ f_*$, where $\rho: \pi_1 (M)\to PSO(1,n)$ is the
defining representation of the hyperbolic manifold $M$, and $f: M \to M$ is an
isometry.  The proof will appear elsewhere.
\pbf {\bf 5.13.} Recall the following classical result of A. Weil [55] and
Garland-Raghunathan [23].
\proclaim{Theorem}  Let $M$ be a complete hyperbolic manifold of finite
volume, $\dim M \ge 3$, and let $\rho: \pi_1(M) \to PSO (1, \w)$ be the
defining representation.  Then $\rho$ is rigid in the following cases:
\item{(a)} $M$ is compact
\item{(b)} $M$ is noncompact and $\dim M \ge 4$.
\endproclaim

Combining this theorem with the argument of 2.2 (called the Vinberg lemma
by Margulis) one get the following important
\proclaim{Corollary} If $M$ is a hyperbolic manifold with $\dim M \ge 3$
and either $M$ is compact or $\dim M \ge 4$, then the defining representation
$\rho$ is defined over a number field.  Alternatively, the lengths of all
closed geodesics of $M$ are algebraic numbers.
\endproclaim

This inspires the following definition.
\pbf{\bf 5.14. Definition.} For $M$ as above, let $F \subset \bbc$ be the
field, generated by the lengths of closed geodesics of $M$.  Then $F$ is
called the field of definition of $M$.  The natural number
$g(M)=[F:\bbq]$
is called the arithmetic genus of $M$.
\pbf{\bf 5.15.} Let $\calh {\cal Y} \calp(n), n\ge 3$, be the set of isometry
classes of  compact hyperbolic manifolds of dimension $n$.  For $n\ge 4$ the
volume function $Vol: \calh {\cal Y }\calp(n) \to \bbr_+$ is proper, i.e.
$\#Vol^{-1} ([0,C]) <\infty$ for any threshold $C > 0$ by the famous theorem
of Wang [54] and Gromov [26] see also [42].  This fails for $n=3$ [52].
However, we state the following conjecture.
\demo{Conjecture} The function $g+Vol: \calh {\cal Y} \calp(3) \to \bbr_+$ is
proper.
Morover, for any number field $F \subset \bbc$ there are but
finitely many uniform lattices in $SL_2(\bbc)$, which are contained in
$SL_2(F)$.
\pbf{\bf 5.16.} Here we show how  to deduce 5.9. (a) from 5.8.  The
conception of the proof borrows a lot from the Gromov's approach to the Mostow
rigidity theorem.  We start with the following lemma.
\pbf{\bf 5.16.1.} \proclaim{Lemma (rigidity of $Vol$)} Let $M$ be a compact
manifold and let
$\rho_t, \; 0\le t \le 1$, be a continuous family of representations of
$\pi_1 (M)$ in $PSO(1,n)$.  Then $Vol (\rho_t)=$ const.
\endproclaim
\demo{Proof} Fix a point $p$ in $\calh^n$ and consider the evaluation map
$v_p: PSO(1,n)\to \calh^n$.  Let $\w$ be the volume form in $\calh^n$.  The
pullback $v_p^* \w$ is a left-invariant form on $PSO(1,n)$ we may consider
the regulator $r (v_p^* \w,\; \;\mu)$ by 3.2.2.  Clearly $r (v_p^* \w,
\mu) =Vol(\mu)\in H^n(M,\bbr)$.

Now, let $\a \in \wedge^n \hbox{\frak\$p\frak\$s\frak\$o}^*(1,n)$ be the
element corresponding
to the left-invariant closed $v_p^* \w$.  Since $\w$ is closed, $\a$ is
a cocycle for the Lie algebra cohomology.  Let $\calp_t$ be the flat principal
$PSO (1,n)$-bundle, corresponding to $\rho_t$.  We may view $\calp_t$ as a
fixed principal bundle with varying flat connection $\w_t \in
\Omega^1(\calp,\hbox{\frak\$p\frak\$s\frak\$o}
(1,n))$.  Then $r (v_p^* \w, \rho_t)$ may be described as the
characteristic class of a $\hbox{\frak\$p\frak\$s\frak\$o}(1,n)$-structure on
$\calp$, corresponding to
the cocycle $\a$,  in the sense of Bott-Halfliger-Bernstein-Rosenfeld [2],
[7].  Since $\hbox{\frak\$p\frak\$s\frak\$o}(1,n)$ is a finite-dimensional
semi-simple algebra, all such
classes are rigid [22].  This proves the lemma.
\pbf {\bf 5.16.2.}  Now let $M$ be a compact hyperbolic manifold, and let
$\rho: \pi_1(M)\to PSO(1,n)$ be the defining representation.  We wish to
prove that $\rho$ is rigid.  Let $\rho_t$ be a path in
$V_{\pi_1(M)}^{PSO(1,n)}$ which starts with $\rho$.  By 5.12.1.,
$Vol(\rho_t)=Vol(\rho)$.  Then by 5.8., $\rho_t$ is conjugate to $\rho\circ
f_{t*}$ for some isometry $f_t$.  But $Iso(M)$ is finite by the Bochner
theorem, so $\rho_t$ is conjugate to $\rho$.
\qed
\pbf{\bf 5.17.}  We will indicate an application of 5.16.1. to the geometry of
representation varieties.  Start with a hyperbolic three-manifold $M$.  Let
$L$ be a link in $M$, homologeous to zero, and let $N\buildrel f \over \to M$
be a $d$-sheet ramified covering along $L$.  Assume that $N$ is itself
hyperbolic (this is often the case by the Thurston theory).  Let $\mu,\nu$
be the defining representations of $\pi_1 (M)$, respectively $\pi_1(N)$ in
$PSL_2 (\bbc)$.  Then by the Gromov inequality, $Vol(N)$ is strictly larger
than $d\cdot Vol(N)$.  Hence, by 5.16, $V^{PSL_2 (\bbc)}_{\pi_1(N)}$ at least
two connected components (in classical topology), containing $\nu$ and $\mu
\circ f_*$, respectively.  Proceeding in this way with $M$ replaced by $N$,
we come to the following
\proclaim{Proposition}  There exists an irreducible compact three-manifold
$P$, such that the representation variety $V^{PSL_2(\bbc)}_{\pi_1(p)}$ has an
arbitrarily large number of components.
\endproclaim
\pbf{\bf 5.18.}  If $M$ is a noncompact  hyperbolic three-manifold, then
5.9(b) fails to be true.  Moreover, the well-known result of Thurston
(see [12]) estimates the dimension of $V_{\pi_1 (M)} (\rho)$ by the number of
cusps in $M$.  We may alter the definition 5.10. as follows: by 5.2.
$V_{\pi_1}(\rho)$ contains a $\bar \bbq$-point $\bar \rho$.  Let $F$ be a
field of the smallest degree such that there exists a $F$-point in
$V_{\pi_1 (M)}(\rho)$; put $g(M)=[F:\bbq]$.  In particular, for an excellent
knot $K\subset S^3$ this defines {\bf the arithmetic genus of the knot $K$}
when applied to the (hyperbolic of finite volume) knot manifold $S^3\backslash
K$.

\noindent\pbf{\bf 5.19} We begin a proof of 5.1. Consider the flat
$\calh^n$-bundle $\calf$ over $M$, corresponding to $\mu$, that is,
$\calf=\tilde M \mathop\times\limits_{\pi_1(M)}
\calh^n$, where $\pi_1(M)$ acts in $\tilde M \times \calh^n \equiv
\calh^n \times \calh^n$ by the diagonal action $(\rho, \mu)$.  Fix a section,
say $s$, of $\calf$.  By the well-known relation between sections and
equivariant maps (see [11], for example), this gives rise to an equivariant
map $\bar s: \tilde M\to \calh^n$ with respect to the actions $\rho$ and
$\mu$.  Let $\Sigma a_i \sigma_i$ be a closed singular chain in $M$.  Let
$\tilde\sigma_i$ be a
lift of $\sigma_i$ to $\tilde M=\calh^n$.  Let $\hat \sigma_i$ be the
Thurston straightening of the singular simplex $\bar s(\tilde \sigma_i)$.
Consider the chain $\Sigma a_i (\tilde \sigma_i, \hat \sigma_i)$ in
$\calh^n\times \calh^n$.  Denote by $p: \calh^n\times \calh^n \to \calf$
the natural projection and consider the chain $b=\Sigma a_i p(\tilde \sigma_i,
\hat \sigma_i)$ in $\calf$.  Let $\pi: \calf\to M$ be the fibration map,
then $\pi(b) =\Sigma a_i \sigma_i$. We claim that $b$ is closed.  This
follows immediately from the description of the straightening process (see
[52]) and the fact that $\Sigma a_i \sigma_i$ is closed in $M$.  Next,
let $\w$ be the volume form in $\calh^n$ and let $\pi_2: \calh^n \times
\calh^n\to\calh^n$ be the projection on the second factor. Then clearly
$Vol(\mu)=\int\limits_{\Sigma_i (\tilde \sigma_i,
\hat \sigma_i)}\pi^*_2 \w \le \Sigma| a_i|\cdot d_n$ where $d_n$ is the
Milnor constant.  Taking  infinum over all chains  representing [M], we
get $Vol(\mu) \le \|[M]\|_g$.
\qed


\noindent {\bf 6. THE CHERN-SIMONS INVARIANT OF ARITHMETIC HYPERBOLIC
THREE-MANIFOLDS}

\noindent{\bf 6.1} In this section we will prove the theorem D.  Recall that
all
arithmetic hyperbolic three-manifolds are constructed as follows.  Let $F$ be
a totally real number field, and let $Q$ be a quadratic form in four
variables, defined over $F$.  Suppose that $Q$ has the signature (1.3) and
that for any nontrivial embedding $\sigma: F \to \bbr$ the form $Q^\sigma$
is negatively defined.  Then $SO(Q)\cap SL_4 (\calo) \subset SO(1,3)$ is
a uniform lattice.
\pbf{\bf 6.2.1.} Fix the identification $PSL_2(\bbc) \approx SO(1,3)$.
It follows that the defining  representation $\rho$ of an arithmetic
hyperbolic three-manifold $M$ is defined over $K=F[\sqrt {-1}]$ where $F$
is totally real.  In particular, all embeddings $\sigma: K\to \bbc$ commute
with the complex conjugation.  Let $\{\sigma_i\}$ be the maximal family
of nonconjugate embeddings.  Then by 4.2., for any representation $\mu:
\pi_1(M)\to SL_2 (\bbc)$
$$
ChS(\mu)= \sum^m_{i=1} \a_i Vol (\sigma_i \circ \mu) \quad  (\hbox{mod} \bbz
\frac{1}{M'})\tag *
$$
Since the defining representation $\rho:\pi_1 (M)\to PSL_2 (\bbc)$ lifts
to $SL_2(\bbc)$ [\ ].  (\*) applies to $\rho$ without any change, so
$$
ChS (\rho)=\sum^m_{i=1} \a_i Vol (\sigma_i \circ \rho)\quad (\hbox{mod} \bbz
\frac{1}{M'})\tag **
$$
Now, applying (*) to the complex-conjugate representation $\bar \rho$ we get
$$
ChS (\bar \rho)=\sum^m_{i=1} \a_i \;Vol(\sii \circ \bar \rho)=\sum^m_{i=1}
\ai \;Vol (\overline{\sii \circ \rho}) (\hbox{mod} \bbz \cdot \frac{1}{M'})
\tag ***
$$,
since $\sii$ commutes with the conjugation.  Now we will use the following
lemma.
\proclaim{{\bf 6.2.2. Lemma}} For any representation $\mu: \pi_1 (M)\to
PSL_2 (\bbc)$ we have
\item{(i)} $ChS(\bar \mu) =ChS(\mu)\;\; (\hbox{mod} \bbz)$
\item{(ii)} $Vol (\bar \mu)=-Vol (\mu)$

\demo{Proof} Let $\w(X,Y,Z)=([X,Y],Z)$ be the canonical 3-form in
$SL_2(\bbc)$.  Then by 3.3, $ChS(\mu)=r(Re \;\w,\; \mu)$ whereas $Vol (\mu)
=r (Im \;\w, \mu)$, and the statement of the lemma follows readily.
\pbf{\bf 6.2.3.}  To finish the proof of the theorem D, we add (**) and
(***) to get $2ChS(\rho)= \sum \ai (Vol (\sii \circ \rho)+Vol (\overline{\sii
\circ \rho}))=0 \;\; (\hbox{mod} \bbz\cdot \frac{1}{M})$, hence $ChS(\rho)\in
\bbq$.
\qed
\pbf{\bf 6.3. Remarks.}  What we actually have proved is $2M'\cdot ChS(\rho)\in
\bbz$, where $M$ is defined by 3.3.  The inspection of 3.3 shows that $M'$ is
just the squared order of the (torsion) second Chern class of the classifying
bundle over $BSL_2 (F[\sqrt{-1}])$.  The latter invariant was intensively
studied in algebraic $K$-theory (see [50], [51] for the references therein
and the connection to the Lichtenbaum conjecture).
\pbf{\bf 6.3.1.}  The arithmetic hyperbolic three-manifolds constitute
relatively ``small'' part of all hyperbolic manifolds; in particular, the set
of volumes
of these manifolds is discrete [10], [5].  However, the corollary 5.13.
shows that in a way {\it all} compact hyperbolic manifolds are
``arithmetic''.  A somewhat deeper look at the proof of the
theorem D above indicates at the level of difficulties in studying the general
case.  If the number field of
definition is not a CM  field, then
the complex conjugation permutes all the embeddings into $\bbc$ in a way we
may not control from the point of view of the arithmetic nature of the
coefficients $\ai$ in (*).  These coefficients depend only on the number
field and are very interesting invariants of it.


\heading REFERENCES \endheading
\item{1.} S.Bauer, Ch.Okonek,{\it The algebraic geometry of representation
spaces
associated to Seifert homology 3-spheres} Math.Ann., {\bf 286} (1990),45--76.
\item{2.}  J. Bernstein, B. Rosenfeld.  {\it Homogeneous spaces of infinite-
dimensional Lie algebras and characteristic classes of foliations}, Sov. Math.
(Uspechi), {\bf 28} (1973),103--138.
\item{3.} S. Bloch, {\it Applications of the dilogarithm function in algebraic
$K$-theory and algebraic geometry}, Proc. Int. Symp. Alg. Geom. Kyoto (1977)
Kinokuniya,
103--114.
\item{4.} A. Borel, {\it Stable real cohomology of arithmetic groups}, Ann.
Sci. Ec. Norm. Super.,(41), {\bf 7} (1974) 235--272.
\item{5.} A. Borel, {\it Commensurability classes and volumes of hyperbolic
3-manifolds}, Ann. Sc. Norm. Super. Pisa, {\bf 8} (1981) 1--33.
\item{6.} A. Borel, J-P. Serre, {\it Cohomologie d'immeubles et de groupes
S-arithm\'etique},
Topology,  {\bf  15} (1976), 211--232.
\item{7.} R. Bott, A. Haefliger, {\it On the characteristic classes of
$\Gamma$-
foliations}, Bull. Amer. Math. Soc. {\bf 78}, 1039--1044.
\item{8.} J. Cheeger and  J. Simons, {\it Differential characters and
geometric invariants}, in ``Geometry and Topology'', Lect. Notes in Math.,
{\bf 1167} (1980) 50--80.
\item{9.} S.-S. Chern and J. Simons,{\it Characteristic forms and geometric
invariants}, Annals Math.,
{\bf 99} (1974), 48--69.
\item{10.} T. Chinburg, {\it Volumes of hyperbolic manifolds}, J. Diff. Geom.,
{\bf 18} (1983) 783--789.
\item{11.} P.Conner, E.Floyd, {\it Differential Periodic Maps},
Springer--Verlag.
\item{12.} M. Culler, P. Shalen, {\it Varieties of group representations and
splittings of 3-manifolds}, Ann. of Math., {\bf 117} (1983) 109--146.
\item{13.} K. Corlette, {\it Flat $G$-bundles with canonical metrics}, J. Diff.
Geom.,
{\bf 28} (1988) 361--382.
\item{14.} K. Corlette, {\it Rigid representations of K\"ahlerian fundamental
groups}, J. Diff. Geom., {\bf 33} (1991) 239--252.
\item{15.} P. Deligne, {\it Th\'eorie de Hodge,I,II,III}, Actes, Congr. Intern.
Math. Nice, 1970, 425--430; Publ.Math.IHES {\bf 40} (1971), 5--58, {\bf 44}
(1974), 5--77.
\item{16.} A. Dolgachev, {\it Automorphic forms and quasihomogeneous
singularities}, Funct. Anal. and Appl. {\bf 9} (1975) 67--68.
\item{17.} S. K. Donaldson, {\it Twisted harmonic maps and the self-duality
equation}, Proc. London Math. Soc., {\bf 55} (1987) 127--131.
\item{18.} J. L. Dupont and C. -H. Sah, {\it Scissors congruences, II}, J.
Pure Appl. Algebra {\bf 25} (1982) 159--195.
\item{19.} J.L. Dupont, {\it The dilogarithms as a characteristic class for
flat bundles}, J.Pure Appl. Algebra, {\bf 44} (1987), 137--164.
\item{20.} H.Esnault, E.Vieweg, {\it Deligne--Beilinson cohomology},
in:Beilinson's Conjecture on Special Values of L-functions,M.Rappoport,
N.Schappacher, P.Schneider, ed.,  Acad.Press, 1988, 43--91.
\item{21.} R. Fintushel and R. Stern, {\it Invariants for homology 3-spheres},
in: ``Geometry of Low-Dimensional Manifolds'', Proc. of the Durham Symposium,
S.K. Donaldson and C.B. Thomas, Ed., Cambridge Univ. Press, 1990, 125--148.
\item{22.} D. B. Fuchs, {\it Cohomology of Infinite-Dimensional Lie Algebras},
Plenum, 1986.
\item{23.} H. Garland and H. Raghunathan, {\it Fundamental domains for
lattices in ($R$)-rank 1 semisimple Lie groups}, Ann. of Math., {\bf 92}
(1970) 279--326.
\item{24.} H. Gillet, C. Soul\'e, {\it Arithmetic Chow groups and differential
characters}, in:``Algebraic $K$-Theory:  connections with Geometry and
Topology'',Kluwer, 1989, 29--68.
\item{25.} W. Goldman, {\it Representations of fundamental groups of surfaces},
in: {\it Geometry and Topology}, ed. J.Alexander and J.Karer
,Lect.Notes in Math. 1167, (1985), 95--117.
\item{26.} M. Gromov, {\it Manifolds of negative curvature}, J. Diff. Geom.,
{\bf 13} (1978) 223--230.
\item{27.} M. Gromov, {\it Hyperbolic manifolds according to Thurston and
Jorgensen}, S\'em. Bourb., 546, Lect. Notes in Math., {\bf 842} (1981) 40--53.
\item{28.} M. Gromov, {\it Hyperbolic groups}, in: {\it Essays in Group Theory}
, ed. S.M.Gersten, Springer--Verlag, 1987, 75--263.
\item{29.} D. Johnson and J. Millson, {\it Deformation spaces associated to
compact hyperbolic manifolds}, in: Discrete Groups in Geometry and Analysis
(papers in honor of G.D. Mostow), Progress in Math., {\bf 67} Birkhauser
(1987) 48-106.
\item{30.} J. D. C. Jones, in preparation.
\item{31.} M. Karoubi, {\it Theorie generale des classes caracteristiques
secondaries}, K-theory, {\bf 4}  (1990/1991)  55--87.
\item{32.} P. Kirk, E. Klassen, {\it Chern-Simons invariants of 3-manifolds
and representation spaces of knot groups}, Math. Ann., {\bf 287} (1990)
343--367.
\item{33.} F. Labourie, {\it Existence d'applications harmonique tordues \`a
valeurs dans les vari\'et\'es \`a courboure n\'egative}, Proc. Amer. Math.
Soc., {\bf 111} (1991) 878--882.
\item{34.} D. MacDuff, {\it Local homology of volume-preserving
diffeomorphisms, I}, Ann. Sci. \'Ecole Norm. Super., {\bf 15} (1982) 609--648.
\item{35.} D. MacDuff,  {\it On groups of volume-preserving diffeomorphisms
and foliations with transverse volume form}, Proc. London Math. Soc., {\bf 43}
(1981) 295--320.
\item{36.} J. Millson, {\it Real vector bundles with discrete structure
group}, Topology, {\bf 18} (1979) 83--80.
\item{37.} J. Milnor, {\it On the 3-dimensional Brieshorn manifolds $M(p,q,
r)$}, in: ``Knots, groups and 3-manifolds'', L.P. Neuwirth, Ed., Princeton
Univ. Press, 1975, 175--225.
\item{38.} N.Mok, {\it Factorization of semisiimple discrete representations
of K\"ahler groups}, Inv. Math. {\bf 110} (1992) 557--614.
\item{39.} S. Morita, {\it Characteristic classes of surface bundles}, Inv.
Math., {\bf 90} (1987) 551--577.
\item{40.} D. Ramakrishnan, {\it Regulators, algebraic cycles and values of
$L$-function},  Contemp. Math., {\bf 83} (1989) 183--310.
\item{41.} A. Reznikov, {\it Harmonic maps, hyperbolic cohomology and higher
Milnor inequalities}, Topology {\bf 32} (1993), 899--907.
\item{42.} A. Reznikov, {\it The volume and the injectivity radius of a
hyperbolic manifold},Topology {\bf 33} (1994), to appear.
\item{43.} A. Reznikov, {\it Symplectic twistor spaces}, Ann. Global Anal. and
Geom.{\bf 11}, No.2(1993), 109--118.
\item{44.} A. Reznikov, {\it Hyperbolic properties of the Chern-Simons
invariant and the topological complexity of homology spheres}, preprint (1993)
\item{45.} A. Reznikov, {\it All regulators of flat bundles are torsion},
preprint
(June,1993)
\item{46.} A. Reznikov, {\it Symplectic Chern--Weil theory},
preprint (1993)
\item{47} J.H. Sampson, {\it Harmonic maps in K\"ahler geometry}, Lect. Notes
in Math. {\bf 1161} (1985) 193--205.
\item{48.} Y.-T. Siu, {\it The complex-analyticity of harmonic maps and the
strong rigidity of compact K\"ahler manifolds}, Ann. of Math. {\bf 112}
(1980) 73--111.
\item{49.} C. Soul\'e, {\it Connexions et classes Caract\'eristiques de
Beilinson}, Contemp. Math., {\bf 83} (1989) 349--376.
\item{50.} C. Soul\'e, {\it Algebraic $K$-theory of the integers}, in:
``Higher Algebraic $K$-theory: an overview'', Lect. Notes in Math. {\bf 1491}
(1992) 100--132.
\item{51.} C.B. Thomas, {\it Characteristic Classes and the Cohomology of
Finite  Groups}, Cambridge Univ. Press, 1986.
\item{52.} W. Thurston, {\it Geometry and Topology of 3-manifolds}, Princeton,
1978.
\item{53.} W. Thurston, unpublished manuscript.
\item{54.} H.C. Wang, {\it Topics in totally discontinuous groups}, in:
Symmetric Spaces, Marcel Dekker,  1972, 460--485.
\item{55.} A. Weil, {\it Discrete subgroups of Lie groups, I, II}, Ann. of
Math. {\bf 72} (1960) 369--384, {\bf  75} (1962) 578--602.
\item{56.} K. Yano,{\it Gromov invariant and $S^1$-action}, J. Fac. Sci. Tokyo,
{\bf 18} (1983) 783--789.


\bye
\end